\documentclass[12pt,oneside,english]{amsart}
\usepackage[T1]{fontenc}
\usepackage[latin9]{inputenc}
\usepackage{geometry}
\usepackage{babel}
\usepackage{units}
\usepackage{mathrsfs}
\usepackage{mathtools}
\usepackage{amsbsy}
\usepackage{amstext}
\usepackage{amsthm}
\usepackage{esint}
\usepackage[unicode=true,pdfusetitle,
 bookmarks=true,bookmarksnumbered=false,bookmarksopen=false,
 breaklinks=false,pdfborder={0 0 1},backref=false,colorlinks=false]
 {hyperref}

\makeatletter
\numberwithin{equation}{section}
\numberwithin{figure}{section}
\theoremstyle{plain}
\newtheorem{conjecture}{\protect\conjecturename}[section]
\theoremstyle{plain}
\newtheorem{thm}{\protect\theoremname}[section]
\theoremstyle{plain}
\newtheorem{lem}{\protect\lemmaname}[section]

\usepackage{babel}

\providecommand{\lemmaname}{Lemma}
\providecommand{\theoremname}{Theorem}

\providecommand{\conjecturename}{Conjecture}

\makeatother

\providecommand{\conjecturename}{Conjecture}
\providecommand{\lemmaname}{Lemma}
\providecommand{\theoremname}{Theorem}

\begin{document}
\title[A conjecture of Baruah and Sarma]{Proof of a conjecture of Baruah and Sarma on sign patterns of certain
infinite products}
\author{Bing He}
\address{School of Mathematics and Statistics, HNP-LAMA, Central South University
\\
 Changsha 410083, Hunan, People's Republic of China}
\email{yuhelingyun@foxmail.com; yuhe123456@foxmail.com}
\author{Xiongze Zhang}
\address{School of Mathematics and Statistics, HNP-LAMA, Central South University
\\
 Changsha 410083, Hunan, People's Republic of China}
\email{17680089436@163.com}
\begin{abstract}
Let 
\[
\sum_{n=0}^{\infty}A(n)q^{n}\coloneqq\frac{(q^{2};q^{5})_{\infty}^{5}(q^{3};q^{5})_{\infty}^{5}}{(q;q^{5})_{\infty}^{5}(q^{4};q^{5})_{\infty}^{5}},
\]

\[
\sum_{n=0}^{\infty}B(n)q^{n}\coloneqq\frac{(q;q^{5})_{\infty}^{5}(q^{4};q^{5})_{\infty}^{5}}{(q^{2};q^{5})_{\infty}^{5}(q^{3};q^{5})_{\infty}^{5}},
\]
and

\[
\sum_{n=0}^{\infty}D(n)q^{n}\coloneqq\frac{(q^{5};q^{25})_{\infty}(q^{20};q^{25})_{\infty}}{(q^{10};q^{25})_{\infty}(q^{15};q^{25})_{\infty}}\frac{(q^{2};q^{5})_{\infty}^{5}(q^{3};q^{5})_{\infty}^{5}}{(q;q^{5})_{\infty}^{5}(q^{4};q^{5})_{\infty}^{5}}
\]
where $(a;q)_{\infty}:=\prod_{k=0}^{\infty}(1-aq^{k})$ and $|q|<1.$
These sequences are closely related to the celebrated Rogers--Ramanujan
continued fraction.

In this paper, we study the sign behavior o of the coefficients $A(n),B(n)$
and $D(n).$ We prove that for all integers $n\geq0,$ 
\begin{align*}
A(5n)<0\quad(n\neq0),\qquad B(5n)<0\quad(n\neq0),\qquad D(5n+1)>0.
\end{align*}
 This confirms a recent conjecture of Baruah and Sarma \cite{BS}.
Our proof is different from the previous method of Baruah and Sarma,
and combines asymptotic coefficient analysis with symbolic computation
for finite case verification.
\end{abstract}

\maketitle

\section{Introduction}

The famous Rogers--Ramanujan identities are given by
\[
\begin{aligned} & \sum_{n=0}^{\infty}\frac{q^{n^{2}}}{(q;q)_{n}}=\frac{1}{(q;q^{5})_{\infty}(q^{4};q^{5})_{\infty}},\\
 & \sum_{n=0}^{\infty}\frac{q^{n^{2}+n}}{(q;q)_{n}}=\frac{1}{(q^{2};q^{5})_{\infty}(q^{3};q^{5})_{\infty}},
\end{aligned}
\]
where, throughout the paper, we use the standard notations for complex
numbers $a$ and $q$ with $|q|<1:$ 
\[
(a;q)_{\infty}\coloneqq\prod_{k=0}^{\infty}(1-aq^{k}),\qquad(a;q)_{n}=\frac{(a;q)_{\infty}}{(aq^{n};q)_{\infty}}.
\]
These identities bridge the theories of partitions, $q$-series, modular
forms, and statistical mechanics. They were first discovered (but
not published) by Rogers \cite{R1} in 1894, independently rediscovered
by Ramanujan before 1913 (cf. Hardy \cite{H}) and later  proved by
Schur \cite{S} in 1917.

The celebrated Rogers--Ramanujan continued fraction is defined as

\[
R(q):=\cfrac{q^{1/5}}{1+\cfrac{q}{1+\cfrac{q^{2}}{1+\cfrac{q^{3}}{1+\ddots}}}},\quad|q|<1.
\]
It was originally discovered by Rogers \cite{R1} in 1894 and later
independently rediscovered by Ramanujan between 1910 and 1913.  Today,
it remains a central object in partitions, modular forms, vertex algebras,
integrable models, and experimental mathematics.

A key connection between the Rogers--Ramanujan identities and the
Rogers-Ramanujan continued fraction is
\[
R(q)=\frac{H(q)}{G(q)}=\frac{(q;q^{5})_{\infty}(q^{4};q^{5})_{\infty}}{(q^{2};q^{5})_{\infty}(q^{3};q^{5})_{\infty}},
\]
where
\begin{align*}
G(q):=\sum_{n=0}^{\infty}\frac{q^{n^{2}}}{(q;q)_{n}},\qquad H(q):=\sum_{n=0}^{\infty}\frac{q^{n^{2}+n}}{(q;q)_{n}}.
\end{align*}
This very important relationship was independently discovered by Rogers
\cite{R1} and Ramanujan \cite{R}.

In 1978, Richmond and Szekeres \cite{RS} studied the Taylor coefficients
of a broad class of infinite products, including $R(q)$ and $1/R(q).$
Define
\[
\frac{1}{R(q)}=\frac{(q^{2};q^{5})_{\infty}(q^{3};q^{5})_{\infty}}{(q;q^{5})_{\infty}(q^{4};q^{5})_{\infty}}:=\sum_{n=0}^{\infty}c(n)q^{n}.
\]
They proved that, as $n\rightarrow\infty,$
\[
c(n)=\frac{\sqrt{2}}{(5n)^{3/4}}\exp\left(\frac{4\pi}{25}\sqrt{5n}\right)\times\left\{ \cos\left(\frac{2\pi}{5}\left(n-\frac{2}{5}\right)\right)+\mathcal{O}(n^{-1/2})\right\} .
\]
This implies that for sufficiently large $n,$
\[
c(5n)>0,\ c(5n+1)>0,\ c(5n+2)<0,\ c(5n+3)<0,\ c(5n+4)<0.
\]

Similarly, define
\[
R(q)=\frac{(q;q^{5})_{\infty}(q^{4};q^{5})_{\infty}}{(q^{2};q^{5})_{\infty}(q^{3};q^{5})_{\infty}}\eqqcolon\sum_{n=0}^{\infty}d(n)q^{n}.
\]
They also showed that as $n\rightarrow\infty,$

\[
d(n)\sim\frac{2^{1/2}}{5^{3/4}}\cos\left(\frac{4\pi}{5}\left(n+\frac{3}{20}\right)\right)n^{-3/4}\exp\left(\frac{4\pi}{5}\sqrt{\frac{n}{5}}\right).
\]
 Thus, for $n$ sufficiently large,

\[
d(5n)>0,\ d(5n+1)<0,\ d(5n+2)>0,\ d(5n+3)<0,\ d(5n+4)<0.
\]

Recently, Baruah and Sarma \cite{BS} investigated the sign behavior
of the coefficients of the higher-order infinite products $R^{5}(q),1/R^{5}(q),R^{5}(q)/R(q^{5})$
and $R(q^{5})/R^{5}(q).$ Let

\[
\frac{1}{R^{5}(q)}=:\sum_{n=0}^{\infty}A(n)q^{n}.
\]
They \cite[Theorem 2]{BS} determined the eventual sign patterns of
$A(n)$ and proved that for all non-negative integers $n,$
\[
A(5n+1)>0,A(5n+2)>0,A(5n+3)>0,A(5n+4)<0.
\]
Let
\begin{align*}
R^{5}(q) & =:\sum_{n=0}^{\infty}B(n)q^{n},\\
\frac{R^{5}(q)}{R(q^{5})} & =:\sum_{n=0}^{\infty}C(n)q^{n},\\
\frac{R(q^{5})}{R^{5}(q)} & =:\sum_{n=0}^{\infty}D(n)q^{n}.
\end{align*}
They \cite[Theorems 3, 4 and 5]{BS} also showed that for all non-negative
integers $n,$ 
\begin{align*}
B(5n+1) & <0,B(5n+2)>0,B(5n+3)<0,B(5n+4)>0,\\
C(5n+1) & <0,C(5n+2)>0,C(5n+3)<0,C(5n+4)>0,C(5n+5)<0,
\end{align*}
and
\[
D(5n+2)>0,D(5n+3)>0,D(5n+4)<0,D(5n+5)<0.
\]
Many interesting congruences satisfied by these coefficients can also
be found in \cite[Theorem 6]{BS}.

At the end of their paper, based on numerical observation, they conjectured
several inequalities on $A(n),B(n)$ and $D(n).$
\begin{conjecture}
\label{c1} \textup{(See \cite[Conjecture 13]{BS})} \footnote{Actually, $A(0)=B(0)=1>0$ and so the two coefficients should be excluded.
Thus, we modify it here.}For all integers $n\geq0,$ we have
\begin{align*}
 & A(5n)<0,\quad n\neq0,\\
 & B(5n)<0,\quad n\neq0,\\
 & D(5n+1)>0.
\end{align*}
\end{conjecture}
Up to now, as far as we are concerned, no proof for this conjecture
has been given.

In this paper, we confirm this conjecture by combining asymptotic
coefficient analysis with finite case verification using symbolic
computation. Our proof is different from the method of Baruah and
Sarma \cite{BS}.
\begin{thm}
\label{t1} Conjecture \ref{c1} is true.
\end{thm}
The paper is organised as follows. In Section \ref{sec:P}, we provide
several auxiliary results that are crucial for the proof of Theorem
\ref{t1}. These include the necessary modular-transformations, integral
estimates and Bessel-function bounds. Section \ref{sec:Pr} supplies
the detailed proof for Theorem \ref{t1}.

\section{\label{sec:P} Preliminaries}

\subsection{Dedekind\textquoteright s eta-function and Jacobi theta function}

The Dedekind eta-function and the Jacobi theta function are defined
as \cite{A,Z}
\[
\eta(\tau)\coloneqq q^{\nicefrac{1}{24}}(q;q)_{\infty},
\]
and 
\[
\vartheta(\varsigma;\tau)\coloneqq\sum_{\nu=\mathbb{Z}+\frac{1}{2}}e^{2\pi i\nu(\varsigma+\frac{1}{2})+\pi i\nu^{2}\tau},
\]
where $q=e^{2\pi i\tau}$ with $\tau\in\mathbb{H}:=\{\tau\in\mathbb{C}|\Im(\tau)>0\}$
and $\varsigma\in\mathbb{C}.$ If $\xi\coloneqq e^{2\pi i\varsigma},$
then, by the Jacobi triple product identity, 
\[
\vartheta(\varsigma;\tau)=-iq^{\nicefrac{1}{8}}\xi^{-\nicefrac{1}{2}}(\xi,\xi^{-1}q;q)_{\infty}.
\]
Let 
\[
\psi(\varsigma;\tau)\coloneqq(\xi,\xi^{-1}q;q)_{\infty}.
\]
Then
\begin{equation}
\psi(\varsigma;\tau)=ie^{-\frac{\pi i\tau}{6}}e^{\pi i\varsigma}\frac{\vartheta(\varsigma;\tau)}{\eta(\tau)}.\label{eq:2-1}
\end{equation}

Let 
\[
\gamma=\left(\begin{array}{cc}
a & b\\
c & d
\end{array}\right)\in\mathrm{SL}_{2}(\mathbb{Z}),\quad c>0.
\]
Recall that the Möbius transformation for $\tau\in\mathbb{H}$ is
defined by
\[
\gamma(\tau)\coloneqq\frac{a\tau+b}{c\tau+d}
\]
For convenience we write 
\[
\gamma^{*}(\tau)\coloneqq\frac{1}{c\tau+d},
\]
and 
\[
\chi(\gamma):=\exp\left(\pi i\left(\frac{a+d}{12c}-s(d,c)-\frac{1}{4}\right)\right),
\]
where $s(d,c)$ is the Dedekind sum defined by 
\[
s(d,c):=\sum_{n(\bmod c)}\left(\left(\frac{dn}{c}\right)\right)\left(\left(\frac{n}{c}\right)\right)
\]
with\footnote{$\left\lfloor x\right\rfloor $ denotes the greatest integer $\leq x.$}

\[
((x)):=\begin{cases}
x-\lfloor x\rfloor-1/2 & \text{if }x\notin\mathbb{Z},\\
0 & \text{if }x\in\mathbb{Z}.
\end{cases}
\]
Then
\begin{equation}
\eta(\gamma(\tau))=\chi(\gamma)(c\tau+d)^{\nicefrac{1}{2}}\eta(\tau)\label{eq:2-2}
\end{equation}
and

\begin{equation}
\vartheta(\varsigma\gamma^{*}(\tau);\gamma(\tau))=\chi(\gamma)^{3}(c\tau+d)^{\nicefrac{1}{2}}e^{\frac{\pi ic\varsigma^{2}}{c\tau+d}}\vartheta(\varsigma;\tau).\label{eq:2-3}
\end{equation}
Furthermore, for integers $A$ and $B,$ we have
\begin{equation}
\vartheta(\varsigma+A\tau+B;\tau)=(-1)^{A+B}e^{-\pi iA^{2}\tau}e^{-2\pi iA\varsigma}\vartheta(\varsigma;\tau).\label{eq:2-4}
\end{equation}

Suppose that $0\leq h<k\leq N$ are three integers such that $\gcd(h,k)=1.$
Let $r<m$ be two positive integers and let $d=\gcd(m,k).$ For convenience,
let $m=dm^{\prime}$ and $k=dk^{\prime}.$ Then $\gcd(m',k')=1.$
This implies that $\gcd(m'h,k')=1.$ Then, there is an integer $\hbar_{m}(h,k)$
such that 
\[
\hbar_{m}(h,k)m'h\equiv-1(\bmod\;k').
\]
Let $b_{m^{\prime}}=\ensuremath{(\hbar_{m}(h,k)m^{\prime}h+1)/k^{\prime}}.$
Then

\[
\gamma_{(m,h,k)}=\begin{pmatrix}\hbar_{m}(h,k) & -b_{m^{\prime}}\\
k^{\prime} & -m^{\prime}h
\end{pmatrix}\in SL_{2}(\mathbb{Z}).
\]

Let 
\[
\tau:=(h+iz)\ensuremath{/}k.
\]
Then
\begin{align*}
\gamma_{(m,h,k)}(m\tau) & =\frac{\hbar_{m}(h,k)\cdot m\frac{h+iz}{dk^{\prime}}-b_{m^{\prime}}}{k^{\prime}\cdot m\frac{h+iz}{dk^{\prime}}-m^{\prime}h}\\
 & =\frac{\hbar_{m}(h,k)m^{\prime}h+i\hbar_{m}(h,k)m^{\prime}z-(\hbar_{m}(h,k)m^{\prime}h+1)}{m^{\prime}hk^{\prime}+ik^{\prime}m^{\prime}z-m^{\prime}hk^{\prime}}\\
 & =\frac{\hbar_{m}(h,k)}{k^{\prime}}+\frac{1}{m^{\prime}k^{\prime}z}i.
\end{align*}
Namely,
\begin{equation}
\gamma_{(m,h,k)}(m\tau)=\frac{\hbar_{m}(h,k)\gcd(m,k)}{k}+\frac{\gcd(m,k)}{mkz}i.\label{eq:2-7}
\end{equation}

Similarly, we get 
\[
\gamma_{(m,h,k)}^{*}(m\tau)=\frac{1}{k^{\prime}\cdot m\frac{h+iz}{dk^{\prime}}-m^{\prime}h}=-\frac{\gcd(m,k)}{mz}i,
\]
so that 
\[
r\tau\gamma_{(m,h,k)}^{*}(m\tau)=\frac{r\gcd(m,k)}{mk}-\frac{rh\gcd(m,k)}{mkz}i.
\]
Therefore, 
\begin{align}
\begin{aligned} & r\tau\gamma_{(m,h,k)}^{*}(m\tau)+\lambda_{m,r}(h,k)\gamma_{(m,h,k)}(m\tau)\\
 & =\frac{r\gcd(m,k)}{mk}+\lambda_{m,r}(h,k)\frac{\hbar_{m}(h,k)\gcd(m,k)}{k}+\lambda_{m,r}^{*}(h,k)\frac{\gcd(m,k)}{mkz}i,
\end{aligned}
\label{eq:2-8}
\end{align}
where\footnote{$\left\lceil x\right\rceil $ denotes the least integer $\geq x.$}
\[
\lambda_{m,r}(h,k):=\left\lceil \frac{rh}{\gcd(m,k)}\right\rceil 
\]
and
\[
\lambda_{m,r}^{*}(h,k):=\lambda_{m,r}(h,k)-\frac{rh}{\gcd(m,k)}.
\]

It follows from \eqref{eq:2-1} that 
\[
\psi(r\tau;m\tau)=ie^{-\frac{\pi im\tau}{6}}e^{\pi ir\tau}\frac{\vartheta(r\tau;m\tau)}{\eta(m\tau)}.
\]
Then, by \eqref{eq:2-2}, \eqref{eq:2-3}, \eqref{eq:2-4} and the
fact $s(-m^{\prime}h,k)=-s(m^{\prime}h,k),$ we have

\begin{align*}
 & \psi(r\tau;m\tau)\\
 & =ie^{-\frac{\pi im\tau}{6}}e^{\pi ir\tau}\chi(\gamma_{(m,h,k)})^{-2}e^{-\frac{\pi ik^{\prime}r^{2}\tau^{2}}{k^{\prime}m\tau-m^{\prime}h}}\frac{\vartheta(r\tau\gamma_{(m,h,k)}^{*}(m\tau);\gamma_{(m,h,k)}(m\tau))}{\eta(\gamma_{(m,h,k)}(m\tau))}\\
 & =i(-1)^{\lambda_{m,r}(h,k)}\frac{\vartheta(r\tau\gamma_{(m,h,k)}^{*}(m\tau)+\lambda_{m,r}(h,k)\gamma_{(m,h,k)}(m\tau);\gamma_{(m,h,k)}(m\tau))}{\eta(\gamma_{(m,h,k)}(m\tau))}\\
 & \quad\times\chi(\gamma_{(m,h,k)})^{-2}e^{-\frac{\pi im\tau}{6}}e^{\pi ir\tau}e^{-\frac{\pi ik^{\prime}r^{2}\tau^{2}}{k^{\prime}m\tau-m^{\prime}h}}e^{\pi i\lambda_{m,r}^{2}(h,k)\gamma_{(m,h,k)}(m\tau)}e^{2\pi i\lambda_{m,r}(h,k)r\tau\gamma_{(m,h,k)}^{*}(m\tau)}\\
 & =i(-1)^{\lambda_{m,r}(h,k)}e^{-2\pi is(-m^{\prime}h,k)}\psi(r\tau\gamma_{(m,h,k)}^{*}(m\tau)+\lambda_{m,r}(h,k)\gamma_{(m,h,k)}(m\tau);\gamma_{(m,h,k)}(m\tau))\\
 & \quad\times\exp\left(\pi i\left(\frac{rh}{k}-\frac{rd}{mk}+\frac{2rd\lambda_{m,r}^{*}(h,k)}{mk}+\frac{\hbar_{m}(h,k)d}{k}\left(\lambda_{m,r}^{2}(h,k)-\lambda_{m,r}(h,k)\right)\right)\right)\\
 & \quad\times\exp\left(\frac{\pi}{12k}\left(\left(2m-12r+\frac{12r^{2}}{m}\right)z-\left(\frac{2d^{2}}{m}+\frac{12d^{2}}{m}\left(\lambda_{m,r}^{*^{2}}(h,k)-\lambda_{m,r}^{*}(h,k)\right)\right)\frac{1}{z}\right)\right)
\end{align*}
Thus,

\begin{align}
\begin{aligned} & \prod_{j=1}^{J}\psi^{\delta_{j}}(r_{j}\tau;m_{j}\tau)\\
 & =i^{\sum_{j=1}^{J}\delta_{j}}(-1)^{\sum_{j=1}^{J}\delta_{j}\lambda_{m_{j},r_{j}}(h,k)}\omega_{h,k}^{2}\varUpsilon_{h,k}\exp\left(\frac{\pi}{12k}\left(\Omega z+\Delta(h,k)z^{-1}\right)\right)\\
 & \times\prod_{j=1}^{J}\psi^{\delta_{j}}\left(r_{j}\tau\gamma_{(m_{j},h,k)}^{*}(m_{j}\tau)+\lambda_{m_{j},r_{j}}(h,k)\gamma_{(m_{j},h,k)}(m_{j}\tau);\gamma_{(m_{j},h,k)}(m_{j}\tau)\right),
\end{aligned}
\label{eq:2-6}
\end{align}
where
\[
\Omega:=\sum_{j=1}^{J}\delta_{j}\left(2m_{j}-12r_{j}+\frac{12r_{j}^{2}}{m_{j}}\right),
\]
\[
\omega_{h,k}:=\exp\left(-\pi i\sum_{j=1}^{J}\delta_{j}\cdot s\left(\frac{m_{j}h}{\gcd(m_{j},k)},\frac{k}{\gcd(m_{j},k)}\right)\right),
\]
\[
\Delta(h,k):=-\sum_{j=1}^{J}\delta_{j}\left(\frac{2\gcd^{2}(m_{j},k)}{m_{j}}+\frac{12\gcd^{2}(m_{j},k)}{m_{j}}\left(\lambda_{m_{j},r_{j}}^{*}(h,k)-\lambda_{m_{j},r_{j}}^{*}(h,k)\right)\right)
\]
and 
\[
\begin{aligned}\varUpsilon_{h,k} & :=\exp\bigg(\pi i\sum_{j=1}^{J}\delta_{j}\bigg(\frac{r_{j}h}{k}-\frac{r_{j}\gcd(m_{j},k)}{m_{j}k}+\frac{2r_{j}\gcd(m_{j},k)\lambda_{m_{j},r_{j}}^{*}(h,k)}{m_{j}k}\\
 & \qquad+\frac{\hbar_{m}(h,k)\gcd(m_{j},k)}{k}\big(\lambda_{m_{j},r_{j}}^{2}(h,k)-\lambda_{m_{j},r_{j}}(h,k)\big)\bigg)\bigg).
\end{aligned}
\]

\subsection{Some inequalities and an integral}

For a Farey fraction $h/k$ of order $N,$ define $\xi_{h,k}\coloneqq[-\theta_{h,k}^{\prime},\theta_{h,k}^{\prime\prime}],$
where $\theta_{h,k}^{\prime},\theta_{h,k}^{\prime\prime}$ are the
distances from $h/k$ to its neighboring mediants. Then
\[
\frac{1}{2kN}\leq\theta_{h,k}^{\prime},\theta_{h,k}^{\prime\prime}\leq\frac{1}{kN}
\]
and so
\[
\frac{1}{kN}\leq|\xi_{h,k}|\leq\frac{2}{kN}.
\]
Also, we have
\[
\Re\left(\frac{1}{z}\right)\geq\Re\left(\frac{1}{k(\rho-i\varphi)}\right)=\frac{1}{k}\Re\left(\frac{\rho+i\varphi}{\rho^{2}+\varphi^{2}}\right)\geq\frac{1}{k}\frac{N^{-2}}{N^{-4}+k^{-2}N^{-2}}\geq\frac{k}{2}.
\]

The modified Bessel function of the first kind is defined as \cite[eq.(4.12.2)]{AAR}
\[
I_{\nu}(x):=\sum_{m=0}^{\infty}\frac{1}{m!\Gamma(m+\nu+1)}\left(\frac{x}{2}\right)^{2m+\nu},
\]
where 
\[
\Gamma(s)=\int_{0}^{\infty}e^{-x}x^{s-1}dx,\quad Re(s)>0.
\]

\begin{lem}
\label{l2.1} For $x\geq3,$ we have 
\[
\frac{1}{10}\frac{e^{x}}{\sqrt{x}}<I_{-1}(x)<\sqrt{\frac{\pi}{8}}\frac{e^{x}}{\sqrt{x}}.
\]
\end{lem}
This inequality is due to Wang \cite[Lemma 2.2]{W}.
\begin{lem}
\label{12.2} Let $(h,k)=1,a\in\mathbb{R}_{>0},b\in\mathbb{R},$ and
define
\[
I\coloneqq\int_{\xi_{h,k}}e^{\frac{\pi}{12k}(bz+az^{-1})}e^{-2\pi in\phi}e^{2\pi n\varrho}d\phi
\]
Then for positive integers $n$ with $n>b/24,$ we have
\end{lem}
\[
I=\frac{2\pi}{k}\left(\frac{24n+b}{a}\right)^{-\frac{1}{2}}I_{-1}\left(\frac{\pi}{6k}\sqrt{a(24n+b)}\right)+E(I)
\]
\textit{where}

\[
\left|E(I)\right|\leq\frac{e^{\frac{\pi a}{3}}e^{2\pi\varrho(n+\frac{b}{24})}}{\pi(n+\frac{b}{24})}.
\]

This is the special $c=0$ case of \cite[Lemma 2.4]{C}.

Let $\eta$ be a positive integer and $\left|\zeta\right|<1,\left|\xi\right|<1,\left|q\right|<1.$
Then

\[
\sum_{n\geq0}\sum_{s\geq0}\sum_{t\geq0}p_{\eta}^{*}(s,t;n)\zeta^{s}\xi^{t}q^{n}=\left(\frac{1}{(\zeta;q)_{\infty}(\xi;q)_{\infty}}\right)^{\eta}.
\]
and
\[
\sum_{n\geq0}\sum_{s\geq0}\sum_{t\boldsymbol{\geq}0}d_{\eta}^{*}(s,t;n)\zeta^{s}\xi^{t}q^{n}=\left((\zeta;q)_{\infty}(\xi;q)_{\infty}\right)^{\eta}.
\]
where $p_{\eta}^{*}(s,t;n)$ is the number of partitions of the positive
integer $n$ into$~\ensuremath{\eta}~$parts, where $s$ parts are
colored red and $t$ parts are colored blue, with $s+t=\eta,$ and
$(-1)^{s+t}d_{\eta}^{*}(s,t;n)$ is the number of partitions of the
positive integer $n$ into$~\ensuremath{\eta}~$distinct parts, with$~s~$of
these parts colored red and $t~$parts colored blue, with $s+t=\eta$.
It is easy to see that 
\[
|d_{\eta}^{*}(s,t;n)|\leq p_{\eta}^{*}(s,t;n)
\]
for all $s,t,n\geq0.$ Then, for $\delta\in\mathbb{Z},$
\[
|(\zeta,\xi;q)_{\infty}^{\delta}|\leq\sum_{n\geq0}\sum_{s\geq0}\sum_{t\geq0}p_{|\delta|}^{*}(s,t;n)|\zeta|^{s}|\xi|^{t}|q|^{n}.
\]
Also, for $0\leq\alpha,\beta,x<1,$

\begin{align}
\begin{aligned} & \sum_{n\geq0}\sum_{s\geq0}\sum_{t\geq0}p_{1}^{*}(s,t;n)\alpha^{s}\beta^{t}x^{n}\\
 & =\frac{1}{(\alpha,\beta;x)_{\infty}}=\exp(-\sum_{k\geq0}\log(1-\alpha x^{k})-\sum_{l\geq0}\log(1-\beta x^{l}))\\
 & \leq\exp\left(\frac{\alpha}{1-\alpha}+\frac{\alpha x}{(1-x)^{2}}+\frac{\beta}{1-\beta}+\frac{\beta x}{(1-x)^{2}}\right),
\end{aligned}
\label{eq: eee}
\end{align}
where we have used the inequalities: for $0\leq x,y<1,$
\[
-\log(1-x)\leq\frac{x}{1-x}
\]
and
\[
-\sum_{k\geq1}\log(1-yx^{k})\leq\frac{xy}{(1-x)^{2}}.
\]

\section{\label{sec:Pr} Proof of Theorem \ref{t1}}

For $f(q)=\sum_{n=0}^{\infty}\alpha(n)q^{n},$ we, by the residue
theorem, have
\[
\alpha(n)=\frac{1}{2\pi i}\ointop_{|q|=r}\frac{f(q)}{q^{n+1}}dq.
\]
Let $\ensuremath{r=e^{-2\pi\rho}}$ with $\ensuremath{\rho=1/N^{2}}$
where $N$ is a sufficiently large positive integer. Set 
\[
q\coloneqq re^{2\pi i\theta}=e^{-2\pi\rho}e^{2\pi i\theta},\quad\theta\coloneqq h/k+\varphi\;\mathrm{on}\;\mathrm{each}\;\xi_{h,k}
\]
 and let $z\coloneqq k(\rho-i\varphi).$ Then
\begin{equation}
\alpha(n)=\sum_{\substack{0\leq h<k\leq N\\
(h,k)=1
}
}e^{-2\pi inh/k}\int_{\xi_{h,k}}f(e^{2\pi i\tau})e^{2\pi n\rho}e^{-2\pi in\varphi}d\varphi\label{eq:2-5}
\end{equation}
In particular, set
\[
f(e^{2\pi i\tau}):=\prod_{j=1}^{J}\psi^{\delta_{j}}(r_{j}\tau;m_{j}\tau).
\]
From \eqref{eq:2-5} and \eqref{eq:2-6} we know that 
\begin{align*}
\alpha(n) & =i^{\sum_{j=1}^{J}\delta_{j}}\sum_{\substack{0\leq h<k\leq N\\
(h,k)=1
}
}e^{-2\pi inh/k}(-1)^{\sum_{j=1}^{J}\delta_{j}\lambda_{m_{j},r_{j}}(h,k)}\omega_{h,k}^{2}\varUpsilon_{h,k}\\
 & \times\int_{\xi_{h,k}}\exp\left(\frac{\pi}{12k}\left(\Omega z+\Delta(h,k)z^{-1}\right)\right)e^{2\pi n\rho}e^{-2\pi in\varphi}\\
 & \times\prod_{j=1}^{J}\psi^{\delta_{j}}\left(r_{j}\tau\gamma_{(m_{j},h,k)}^{*}(m_{j}\tau)+\lambda_{m_{j},r_{j}}(h,k)\gamma_{(m_{j},h,k)}(m_{j}\tau);\gamma_{(m_{j},h,k)}(m_{j}\tau)\right)d\varphi
\end{align*}
Let $L=\mathrm{lcm}(m_{1},\ldots,m_{j})$ and let $1\leq l\leq L,0\leq\aleph<l$
such that $k\equiv l\ (\bmod\;L)$ and $h\equiv\aleph\ (\bmod\ l).$
Then for all $j=1,2,\ldots,J,$

\[
\gcd(m_{j},k)=\gcd(m_{j},l),
\]
and 
\begin{equation}
\left(\lambda_{m_{j},r_{j}}^{*}(h,k)\right)^{2}=\left(\lambda_{m_{j},r_{j}}^{*}(\aleph,l)\right)^{2}.\label{eq:3-1}
\end{equation}
These imply that 
\[
\Delta(h,k)=\Delta(\aleph,l).
\]
Thus,

\begin{align*}
\alpha(n) & =i^{\sum_{j=1}^{J}\delta_{j}}\mathop{\sum_{1\leq l\leq L}\sum_{0\leq\aleph<l}\sum_{\substack{1\leq k\leq N\\
k\equiv l(\bmod L)
}
}\sum_{\substack{0\leq h<k\\
(h,k)=1\\
h\equiv\aleph(\bmod l)
}
}}e^{-2\pi inh/k}\\
 & \times(-1)^{\sum_{j=1}^{J}\delta_{j}\lambda_{m_{j},r_{j}}(h,k)}\omega_{h,k}^{2}\varUpsilon_{h,k}\times\int_{\xi_{h,k}}\exp\left(\frac{\pi}{12k}\left(\Omega z+\Delta(\aleph,l)z^{-1}\right)\right)e^{2\pi n\rho}e^{-2\pi in\varphi}\\
 & \times\prod_{j=1}^{J}\psi^{\delta_{j}}(r_{j}\tau\gamma_{(m_{j},h,k)}^{*}(m_{j}\tau)+\lambda_{m_{j},r_{j}}(h,k)\gamma_{(m_{j},h,k)}(m_{j}\tau);\gamma_{(m_{j},h,k)}(m_{j}\tau))d\varphi\\
 & \eqqcolon i^{\sum_{j=1}^{J}\delta_{j}}\sum_{1\leq l\leq L}\sum_{0\leq\aleph<l}S_{\aleph,l}.
\end{align*}

Let 
\[
\ensuremath{\tilde{\varsigma_{j}}=r_{j}\tau\gamma_{(m_{j},h,k)}^{*}(m_{j}\tau)+\lambda_{m_{j},r_{j}}(h,k)\gamma_{(m_{j},h,k)}(m_{j}\tau)}
\]
and
\[
\ensuremath{\tilde{\tau_{j}}=\gamma_{(m_{j},h,k)}(m_{j}\tau)}.
\]
Then, by \eqref{eq:2-7} and \eqref{eq:2-8},

\begin{align}
\tilde{\varsigma_{j}} & =\frac{r_{j}\gcd(mj,k)}{m_{j}k}+\lambda_{m_{j},r_{j}}(h,k)\frac{\hbar_{m_{j}}(h,k)\gcd(m_{j},k)}{k}+\lambda_{m_{j},r_{j}}^{*}(h,k)\frac{\gcd(m_{j},k)}{m_{j}kz}i\label{eq:3-2}
\end{align}
and
\[
\tilde{\tau_{j}}=\frac{\hbar_{m_{j}}(h,k)\gcd(m_{j},k)}{k}+\frac{\gcd(m_{j},k)}{m_{j}kz}i.
\]

For fixed $\aleph$ and $l$ with $1\leq l\leq L$ and $0\leq\aleph<l,$
let
\[
J_{\aleph,l}^{*}=\left\{ j^{*}|1\leq j^{*}\leq J,\lambda_{m_{j^{*}},r_{j^{*}}}^{*}(\aleph,l)=0\right\} 
\]
and
\[
J_{\aleph,l}^{**}=\left\{ j^{**}|1\leq j^{**}\leq J,\lambda_{m_{j^{**}},r_{j^{**}}}^{*}(\aleph,l)\neq0\right\} .
\]
It follows from \eqref{eq:3-1} that the indices set $\left\{ j=1,2,\ldots,J\right\} $
can be splited into these two disjoint sets:
\[
\left\{ j=1,2,\ldots,J\right\} =J_{\aleph,l}^{*}\cup J_{\aleph,l}^{**}.
\]

\subsection{On $\psi^{\delta}(\varsigma;\tau)$}

For $j^{*}\in J_{\aleph,l}^{*}$ we have
\[
\lambda_{m_{j^{*}},r_{j^{*}}}^{*}(\aleph,l)=0.
\]
Then, by \eqref{eq:3-2},
\[
\tilde{\varsigma_{j^{*}}}=\frac{r\gcd(mj^{*},k)}{m_{j^{*}}k}+\lambda_{m_{j^{*}},r_{j^{*}}}(h,k)\frac{\hbar_{m_{j^{*}}}(h,k)\gcd(m_{j^{*}},k)}{k}.
\]
We now write $m=m_{j^{*}},r=r_{j^{*}},$ and $d=\gcd(m,k),m=dm^{\prime},k=dk^{\prime}.$
Then, we have

\begin{align*}
\tilde{\varsigma_{j}} & =\frac{rd}{mk}+\lambda_{m,r}(h,k)\frac{\hbar_{m}(h,k)d}{k}\\
 & =\frac{rd}{mk}+\frac{rh}{d}\frac{\hbar_{m}(h,k)d}{k}\\
 & =\frac{r(1+\hbar_{m}(h,k)m^{\prime}h)}{m^{\prime}k}=\frac{b_{m^{\prime}}r}{dm^{\prime}}.
\end{align*}

Suppose that $\tilde{\varsigma_{j}}$ is an integer, notice that $\lambda_{m,r}^{*}(\aleph,l)=\lambda_{m,r}^{*}(h,k)=0.$
Then $d$ divides $rh,$ so that $d|r,$ since $(h,k)=1,k=dk^{\prime}.$
Since

\[
b_{m^{\prime}}=\frac{1+\hbar_{m}(h,k)m^{\prime}h}{k^{\prime}},
\]
we have $\gcd(b_{m^{\prime}},m^{\prime})=1.$ This implies that $m^{\prime}|\frac{r}{d},$
so that $m|r.$ This contradicts $1\leq r\leq m-1.$ Hence $\tilde{\varsigma_{j}}$
is not an integer.

 It is easy to see that the choice of $\hbar_{m}(h,k)$ does not
affect the values of $\exp(2\pi i(\tilde{\varsigma_{j^{*}}})),$
$\Pi_{h,k}$ and $\varUpsilon_{h,k}.$ Observe that $0\leq\Im(\tilde{\varsigma_{j}})<\Im(\tilde{\tau_{j}}).$
We write

\begin{align}
\prod_{j=1}^{J}\psi^{\delta_{j}}(\tilde{\varsigma_{j}};\tilde{\tau_{j}}) & =\prod_{j^{*}\in J^{*}}(1-e^{2\pi i\tilde{\varsigma_{j^{*}}}})^{\delta_{j^{*}}}\times\prod_{j^{*}\in J^{*}}(e^{2\pi i\tilde{(\varsigma_{j^{*}}}+\tilde{\tau_{j^{*}}})},e^{2\pi i(\tilde{\tau_{j^{*}}}-\tilde{\varsigma_{j^{*}}})};e^{2\pi i\tilde{\tau_{j^{*}}}})_{\infty}^{\delta_{j^{*}}}\label{eq: 3.10}\\
 & \times\prod_{j^{**}\in J^{**}}(e^{2\pi i\tilde{\varsigma_{j^{**}}}},e^{2\pi i(\tilde{\tau_{j^{**}}}-\tilde{\varsigma_{j^{**}}})};e^{2\pi i\tilde{\tau_{j^{**}}}})_{\infty}^{\delta_{j^{**}}}\nonumber 
\end{align}

\subsection{Proof of $A(5n)<0$}

Let 
\[
\mathcal{L}_{>0}:=\{(\aleph,\ell):1\leq\ell\leq L,\ 0\leq\aleph<\ell,\ \Delta(\aleph,\ell)>0\}.
\]
Note that 
\[
\frac{1}{R(q)^{5}}=\psi(2\tau;5\tau)^{5}\psi(\tau;5\tau)^{-5}.
\]
 Then 
\[
m_{1}=m_{2}=5,~r_{1}=2~r_{2}=1,\delta_{1}=5~\delta_{2}=-5
\]
 Hence $L=5$ and $\Omega=-24.$ We compute that $\mathcal{L}_{>0}=\left\{ (1,5),(4,5)\right\} .$

\subsubsection{The case $\Delta(\aleph,l)\protect\leq0$}

Let

\[
S_{\aleph,l}^{(1)}:=\begin{cases}
S_{\aleph,l}, & \mathrm{if}\quad\Delta(\aleph,l)\leq0,\\
\text{0}, & \mathrm{if}\quad\Delta(\aleph,l)>0,
\end{cases}
\]
and
\[
\alpha^{(1)}(n)=i^{\sum_{j=1}^{J}\delta_{j}}\sum_{1\leq l\leq L}\sum_{0\leq\aleph<l}S_{\aleph,l}^{(1)}.
\]
 Notice that when $1\leq l<5,$
\[
\lambda_{m_{j},r_{j}}^{*}(\aleph,l)=0,(m_{j},l)=1\quad\mathrm{for}\quad j=1,2).
\]
 Since
\[
\prod_{j=1}^{2}\psi^{\delta_{j}}(\tilde{\varsigma_{j}};\tilde{\tau_{j}})=\prod_{j=1}^{2}(1-e^{2\pi i\tilde{\varsigma_{j}}})^{\delta_{j}}\times\prod_{j=1}^{2}(e^{2\pi i\tilde{(\varsigma_{j}}+\tilde{\tau_{j}})},e^{2\pi i(\tilde{\tau_{j}}-\tilde{\varsigma_{j}})};e^{2\pi i\tilde{\tau_{j}}})_{\infty}^{\delta_{j}},
\]
we have
\[
\tilde{\varsigma_{j}}=\frac{b_{5}r_{j}}{m_{j}}=\frac{b_{5}r_{j}}{5},\quad j=1,2.
\]
 Since $\frac{b_{5}r_{j}}{5}$ is not an integer, we can easily get
\begin{equation}
|1-e^{\frac{2\pi i}{5}}|\leq|1-e^{2\pi i\frac{b_{5}r_{j}}{5}}|\leq2,\qquad j=1,2.\label{eq3.1}
\end{equation}

It follows from \eqref{eq: eee} that, for $j=1,2,$ 
\begin{align*}
 & \left|(e^{2\pi i\tilde{(\varsigma_{j}}+\tilde{\tau_{j}})},e^{2\pi i(\tilde{\tau_{j}}-\tilde{\varsigma_{j}})};e^{2\pi i\tilde{\tau_{j}}})_{\infty}^{\delta_{j}}\right|\\
 & \leq\sum_{n\geq0}\sum_{s\geq0}\sum_{t\geq0}p_{5}^{*}(s,t;n)e^{-2\pi\Im(\tilde{\tau_{j}})(s+t+n)}\\
 & =\left(\frac{1}{(e^{-2\pi\Im(\tilde{\tau_{j}})},e^{-2\pi\Im(\tilde{\tau_{j}})};e^{-2\pi\Im(\tilde{\tau_{j}})})_{\infty}}\right)^{5}\\
 & \leq\exp\left[10\left(\frac{e^{-2\pi\Im(\tilde{\tau_{j}})}}{1-e^{-2\pi\Im(\tilde{\tau_{j}})}}+\frac{e^{-4\pi\Im(\tilde{\tau_{j}})}}{(1-e^{-2\pi\Im(\tilde{\tau_{j}})})^{2}}\right)\right]\\
 & =\exp\left[\frac{10e^{-2\pi\Im(\tilde{\tau_{j}})}}{(1-e^{-2\pi\Im(\tilde{\tau_{j}})})^{2}}\right].
\end{align*}
Since

\[
\Im(\tilde{\tau_{j}})=\frac{\gcd(5,l)}{5k}\Re\left(\frac{1}{z}\right)\geq\frac{1}{10},\quad j=1,2,
\]
 we have
\begin{equation}
\exp\left[\frac{10e^{-2\pi\Im(\tilde{\tau_{j}})}}{(1-e^{-2\pi\Im(\tilde{\tau_{j}})})^{2}}\right]\leq\exp\frac{10e^{-\frac{\pi}{5}}}{(1-e^{-\frac{\pi}{5}})^{2}}<e^{25},\quad j=1,2.\label{eq3.2}
\end{equation}

When $\Delta(\aleph,l)\leq0$, we have $\Re(z)=\frac{k}{N^{2}}$ and
$\Re(z^{-1})\geq\frac{k}{2}.$ Then
\begin{align*}
\begin{aligned} & \exp\left(\frac{\pi}{12k}\left(\Omega\Re(z)+\Delta(\aleph,\ell)\Re(z^{-1})\right)\right)\\
 & \leq\exp\left(\frac{\pi}{12k}\left(\Omega\frac{k}{N^{2}}+\Delta(\aleph,\ell)\frac{k}{2}\right)\right)\\
 & =\exp\left(\frac{\pi\rho\Omega}{12N^{2}}\right)\exp\left(\frac{\pi\Delta(\aleph,\ell)}{24}\right).
\end{aligned}
\end{align*}
This, tother with \eqref{eq3.1} and \eqref{eq3.2}, implies that

\begin{align}
\begin{aligned} & \left|i^{\sum_{j=1}^{J}\delta_{j}}\sum_{1\leq l\leq4}\sum_{0\leq\aleph<l}S_{\aleph,l}^{(1)}\right|\\
 & \leq\sum_{1\leq l\leq4}\sum_{0\leq\aleph<l}\left|S_{\aleph,l}\right|\\
 & =\mathop{\sum_{1\leq l\leq4}\sum_{0\leq\aleph<l}\sum_{\substack{1\leq k\leq N\\
k\equiv l(\bmod5)
}
}\sum_{\substack{0\leq h<k\\
\gcd(h,k)=1\\
h\equiv\aleph(\bmod l)
}
}}\left|\int_{\xi_{h,k}}\exp\left(\frac{\pi}{12k}\left(\Omega z+\Delta(\aleph,l)z^{-1}\right)\right)\right.\\
 & \quad\times\left.\prod_{j=1}^{2}\psi^{\delta_{j}}(\tilde{\varsigma}_{j};\tilde{\tau}_{j})e^{2\pi n\rho}e^{-2\pi in\varphi}d\varphi\right|\\
 & <e^{50}\times2{}^{5}\times(|1-e^{\frac{2\pi i}{5}}|)^{-5}\\
 & \quad\times\mathop{\sum_{\substack{1\leq l\leq4\\
0\leq\aleph<l
}
}\sum_{\substack{1\leq k\leq N\\
k\equiv l(\bmod5)
}
}\sum_{\substack{0\leq h<k\\
\gcd(h,k)=1\\
h\equiv\aleph(\bmod l)
}
}}\left|\int_{\xi_{h,k}}\exp\left(\frac{\pi}{12k}\left(\Omega z+\Delta(\aleph,l)z^{-1}\right)\right)e^{2\pi n\rho}e^{-2\pi in\varphi}d\varphi\right|\\
 & <e^{54}\sum_{\substack{1\leq l\leq4\\
0\leq\aleph<l
}
}\sum_{\substack{1\leq k\leq N\\
k\equiv l(\bmod5)
}
}\sum_{\substack{0\leq h<k\\
\gcd(h,k)=1\\
h\equiv\aleph(\bmod l)
}
}e^{2\pi\rho(n-1)}\frac{2}{KN}.
\end{aligned}
\label{eq:3.22}
\end{align}

When $l=5$, since 
\[
\Delta(2,5)=\Delta(3,5)=-24,\lambda_{m_{j},r_{j}}^{*}(\mathscr{\aleph},5)\neq0,~(m_{j},5)=5,\quad j=1,2,
\]
we have
\[
\left|\prod_{j=1}^{2}\psi^{\delta_{j}}(\tilde{\varsigma_{j}};\tilde{\tau_{j}})\right|=\prod_{j=1}^{2}\left|(e^{2\pi i\tilde{\varsigma_{j}}},e^{2\pi i(\tilde{\tau_{j}}-\tilde{\varsigma_{j}})};e^{2\pi i\tilde{\tau_{j}}})_{\infty}^{\delta_{j}}\right|.
\]
Note that for $j=1,2,$
\[
\lambda_{m_{j},r_{j}}^{*}(\aleph,l)\geq\frac{1}{5},\quad\Im(\tilde{\varsigma_{j}})=\lambda_{m_{j},r_{j}}^{*}(\aleph,l)\frac{5}{k}\Re\left(\frac{1}{z}\right)\geq\frac{5}{2}\lambda_{m_{j},r_{j}}^{*}(\aleph,l)\geq\frac{1}{2}
\]
and 
\[
\Im(\tilde{\tau_{j}})=\frac{5}{k}\Re\left(\frac{1}{z}\right)\geq\frac{5}{2}.
\]
Then, for $j=1,2$, we have
\[
(e^{-2\pi\Im(\tilde{\tau_{j}})})^{\lambda_{m_{j},r_{j}}^{*}(\aleph,l)}=e^{-2\pi\frac{5}{k}\Re(\frac{1}{z})\times\lambda_{m_{j},r_{j}}^{*}(\aleph,l)}=(e^{-2\pi\frac{5}{k}\Re(\frac{1}{z})\lambda_{m_{j},r_{j}}^{*}(\aleph,l)})=e^{-2\pi\Im(\tilde{\varsigma_{j}})}
\]
and
\begin{align}
\begin{aligned}2e^{-\pi\Im(\tilde{\tau_{j}})} & \leq e^{-2\pi\Im(\tilde{\varsigma_{j}})}+e^{-2\pi\Im(\tilde{\tau_{j}}-\tilde{\varsigma_{j}})}\\
 & =(e^{-2\pi\Im(\tilde{\tau_{j}})})^{\lambda_{m_{j},r_{j}}^{*}(\aleph,l)}+\frac{e^{-2\pi\Im(\tilde{\tau_{j}})}}{(e^{-2\pi\Im(\tilde{\tau_{j}})})^{\lambda_{m_{j},r_{j}}^{*}(\aleph,l)}}\\
 & \leq e^{\frac{-2\pi\Im(\tilde{\tau_{j})}}{5}}+e^{\frac{-8\pi\Im(\tilde{\tau_{j})}}{5}}\leq e^{-\pi}+e^{-4\pi}.
\end{aligned}
\label{eq: 3.7}
\end{align}
It follows that 
\begin{align}
\begin{aligned} & \left|(e^{2\pi i\tilde{\varsigma_{j}}},e^{2\pi i(\tilde{\tau_{j}}-\tilde{\varsigma_{j}})};e^{2\pi i\tilde{\tau_{j}}})_{\infty}^{\delta_{j}}\right|\\
 & \leq\sum_{n\geq0}\sum_{s\geq0}\sum_{t\geq0}p_{5}^{*}(s,t;n)e^{-2\pi\Im(\tilde{\varsigma_{j}})s}e^{-2\pi\Im(\tilde{\tau_{j}}-\tilde{\varsigma_{j}})t}e^{-2\pi\Im(\tilde{\tau_{j}})n}\\
 & \leq\exp\left[5\left(\frac{e^{-2\pi\Im(\tilde{\varsigma_{j}})}}{1-e^{-2\pi\Im(\tilde{\varsigma_{j}})}}+\frac{e^{-2\pi\Im(\tilde{\varsigma_{j}})}e^{-2\pi\Im(\tilde{\tau_{j}})}}{(1-e^{-2\pi\Im(\tilde{\tau_{j}})})^{2}}\right.\right.\\
 & \quad\left.\left.+\frac{e^{-2\pi\Im(\tilde{\tau_{j}}-\tilde{\varsigma_{j}})}}{1-e^{-2\pi\Im(\tilde{\tau_{j}}-\tilde{\varsigma_{j}})}}+\frac{e^{-2\pi\Im(\tilde{\tau_{j}}-\tilde{\varsigma_{j}})}e^{-2\pi\Im(\tilde{\tau_{j}})}}{(1-e^{-2\pi\Im(\tilde{\tau_{j}})})^{2}}\right)\right]\\
 & =\exp\left[5\left(\frac{(e^{-2\pi\Im(\tilde{\varsigma_{j}})}+e^{-2\pi\Im(\tilde{\tau_{j}}-\tilde{\varsigma_{j}})})e^{-2\pi\Im(\tilde{\tau_{j}})}}{(1-e^{-2\pi\Im(\tilde{\tau_{j}})})^{2}}\right.\right.\\
 & \quad\left.\left.+\frac{(e^{-2\pi\Im(\tilde{\varsigma_{j}})}+e^{-2\pi\Im(\tilde{\tau_{j}}-\tilde{\varsigma_{j}})})-2e^{-2\pi\Im(\tilde{\tau_{j}})}}{1-(e^{-2\pi\Im(\tilde{\varsigma_{j}})}+e^{-2\pi\Im(\tilde{\tau_{j}}-\tilde{\varsigma_{j}})})+e^{-2\pi\Im(\tilde{\tau_{j}})}}\right)\right].
\end{aligned}
\label{eq: 3.6}
\end{align}
Combining \eqref{eq: 3.6} and \eqref{eq: 3.7}, we get
\begin{equation}
\begin{aligned} & \left|(e^{2\pi i\tilde{\varsigma_{j}}},e^{2\pi i(\tilde{\tau_{j}}-\tilde{\varsigma_{j}})};e^{2\pi i\tilde{\tau_{j}}})_{\infty}^{\delta_{j}}\right|\\
 & \leq\exp\left[5\left(\frac{(e^{-\pi}+e^{-4\pi})e^{-5\pi}}{(1-e^{-5\pi})^{2}}+\frac{(e^{-\pi}+e^{-4\pi})-2e^{-5\pi}}{1-(e^{-\pi}+e^{-4\pi})+e^{-5\pi}}\right)\right]\\
 & <2
\end{aligned}
\label{eq: 3.12}
\end{equation}
for $j=1,2.$ This implies that 
\[
\left|\prod_{j=1}^{2}\psi^{\delta_{j}}(\tilde{\varsigma_{j}};\tilde{\tau_{j}})\right|<2{}^{2}=4.
\]
From this we obtain 
\begin{align}
|S_{2,5}|+|S_{3,5}| & \leq4\sum_{\substack{1\leq k\leq N\\
k\equiv0(\bmod5)
}
}\sum_{\substack{0\leq h<k\\
\gcd(h,k)=1\\
h\equiv2,3(\bmod5)
}
}e^{2\pi\rho(n-1)}\frac{2}{KN}.\label{eq:3.23}
\end{align}
Thus, by \eqref{eq:3.22} and \eqref{eq:3.23}, 
\[
|\alpha^{(1)}(n)|<e^{54}\sum_{\substack{0\leq h<k\leq N\\
(h,k)=1
}
}e^{2\pi\rho(n-1)}\frac{2}{KN}<2e^{2\pi\rho(n-1)+54}.
\]

\subsubsection{The case $\Delta(\aleph,l)>0$}

Let
\[
S_{\aleph,l}^{(2)}:=\begin{cases}
0, & \Delta(\aleph,l)\leq0,\\
\text{\ensuremath{S_{\aleph,l}}}, & \Delta(\aleph,l)>0.
\end{cases}
\]
Now we split $S_{\aleph,l}^{(2)}$ into two parts $A_{1}$ and $A_{2}:$
\[
S_{\aleph,l}^{(2)}=A_{1}+A_{2},
\]
 where

\begin{align*}
A_{1} & \coloneqq\sum_{\substack{1\leq k\leq N\\
k\equiv0(\bmod5)
}
}\sum_{\substack{0\leq h<k\\
\gcd(h,k)=1\\
h\equiv1,4(\bmod5)
}
}e^{-2\pi inh/k}(-1)^{\sum_{j=1}^{J}\delta_{j}\lambda_{m_{j},r_{j}}(h,k)}\omega_{h,k}^{2}\varUpsilon_{h,k}\\
 & \quad\times\int_{\xi_{h,k}}\exp(\frac{\pi}{12k}(\Omega z+\Delta(h,k)z^{-1}))\Pi_{h,k}e^{-2\pi in\phi}e^{2\pi n\varrho}d\phi
\end{align*}
and

\begin{align*}
A_{2} & \coloneqq\sum_{\substack{1\leq k\leq N\\
k\equiv0(\bmod5)
}
}\sum_{\substack{0\leq h<k\\
\gcd(h,k)=1\\
h\equiv1,4(\bmod5)
}
}e^{-2\pi inh/k}(-1)^{\sum_{j=1}^{J}\delta_{j}\lambda_{m_{j},r_{j}}(h,k)}\omega_{h,k}^{2}\varUpsilon_{h,k}\\
 & \quad\times\int_{\xi_{h,k}}\exp(\frac{\pi}{12k}(\Omega z+\Delta(h,k)z^{-1}))(\prod_{j=1}^{2}\psi^{\delta_{j}}(\tilde{\varsigma_{j}};\tilde{\tau_{j}})-\Pi_{h,k})e^{-2\pi in\phi}e^{2\pi n\varrho}d\phi
\end{align*}
with 
\[
\Pi_{h,k}:=\begin{cases}
\prod_{j\in K^{*}}\bigg(1-\exp\left(2\pi i\frac{r_{j}\gcd(m_{j},k)+r_{j}\hbar_{m_{j}}(h,k)m_{j}h}{m_{j}k}\right)\bigg)^{\delta_{j}}, & \mathrm{if}\;\lambda_{m_{j},r_{j}}^{*}(h,k)=0,\\
1, & \mathrm{otherwise},
\end{cases}
\]
and
\[
K^{*}=\left\{ j|1\leq j\leq J,\;\lambda_{m_{j},r_{j}}^{*}(h,k)=0\right\} .
\]

Notice that 
\[
\Pi_{h,k}=1,\Delta(1,5)=\Delta(4,5)=24.
\]
Then

\begin{align*}
|A_{2}| & \leq\sum_{\substack{1\leq k\leq N\\
k\equiv0(\bmod5)
}
}\sum_{\substack{0\leq h<k\\
(h,k)=1\\
h\equiv1,4(\bmod5)
}
}\int_{\xi_{h,k}}\left|\exp\left(\frac{\pi}{12k}\left(-24z+24z^{-1}\right)\right)\right.\\
 & \left.\times\left(\prod_{j=1}^{2}\psi^{\delta_{j}}(\tilde{\varsigma}_{j};\tilde{\tau}_{j})-1\right)e^{2\pi n\varrho}\right|d\phi\\
 & \leq\sum_{\substack{1\leq k\leq N\\
k\equiv0(\bmod5)
}
}\sum_{\substack{0\leq h<k\\
\gcd(h,k)=1\\
h\equiv1,4(\bmod5)
}
}e^{2\pi\varrho(n-1)}\int_{\xi_{h,k}}\exp\left(\frac{2\pi}{k}\Re\left(\frac{1}{z}\right)\right)\left|\prod_{j=1}^{2}\psi^{\delta_{j}}(\tilde{\varsigma}_{j};\tilde{\tau}_{j})-1\right|d\phi
\end{align*}

It follows from \eqref{eq: 3.10} that for $j=1,2,$
\begin{align*}
\begin{aligned} & \left|\prod_{j=1}^{2}\psi^{\delta_{j}}(\tilde{\varsigma}_{j};\tilde{\tau}_{j})-1\right|\\
 & =\prod_{j=1}^{2}\left|\left(e^{2\pi i\tilde{\varsigma}_{j}},e^{2\pi i(\tilde{\tau}_{j}-\tilde{\varsigma}_{j})};e^{2\pi i\tilde{\tau}_{j}}\right)_{\infty}^{\delta_{j}}-1\right|\\
 & \leq\sum_{\substack{n,s,t\in\mathbb{Z}_{\geq0}^{2}\\
(n,s,t)\neq0
}
}\prod_{j=1}^{2}p_{5}^{*}\left(s,t;n\right)\left|e^{-2\pi\Im\left(\tilde{\varsigma_{j}}\right)s_{j}}\right|\left|e^{-2\pi\Im\left(\tilde{\tau_{j}}-\tilde{\varsigma_{j}}\right)t_{j}}\right|\left|e^{-2\pi\Im\left(\tilde{\tau_{j}}\right)n_{j}}\right|\\
 & =\sum_{\substack{n,s,t\in\mathbb{Z}_{\geq0}^{2}\\
(n,s,t)\neq0
}
}\left(\prod_{j=1}^{2}p_{5}^{*}(s,t;n)\right)\exp\left(-2\pi\frac{\Re(z^{-1})}{k}\right.\\
 & \quad\left.\times5\sum_{j=1}^{2}\left(\lambda_{m_{j},r_{j}}^{*}(\aleph,5)s_{j}+\left(1-\lambda_{m_{j},r_{j}}^{*}(\aleph,5)\right)t_{j}+n_{j}\right)\right).
\end{aligned}
\end{align*}
Hence, for $j=1,2$
\begin{align*}
 & \exp\left(\frac{2\pi}{k}\Re\left(\frac{1}{z}\right)\right)\left|\prod_{j=1}^{2}\psi^{\delta_{j}}(\tilde{\varsigma_{j}};\tilde{\tau_{j}})-1\right|\\
 & \leq\sum_{\substack{n,s,t\in\mathbb{Z}_{\geq0}^{2}\\
(n,s,t)\neq0
}
}\left(\prod_{j=1}^{2}p_{5}^{*}(s,t;n)\right)\exp\left[-2\pi\frac{\Re(z^{-1})}{k}\right.\\
 & \left.\times\left(5\sum_{j=1}^{2}\left(\lambda_{m_{j},r_{j}}^{*}(\aleph,5)s_{j}+\left(1-\lambda_{m_{j},r_{j}}^{*}(\aleph,5)\right)t_{j}+n_{j}\right)-1\right)\right].
\end{align*}

Since 
\[
5\sum_{j=1}^{2}\left(\lambda_{m_{j},r_{j}}^{*}(\aleph,5)s_{j}+(1-\lambda_{m_{j},r_{j}}^{*}(\aleph,5))t_{j}+n_{j}\right)\geq1
\]
we know that

\[
\exp\left(\frac{2\pi}{k}\Re\left(\frac{1}{z}\right)\right)\left|\prod_{j=1}^{2}\psi^{\delta_{j}}(\tilde{\varsigma_{j}};\tilde{\tau_{j}})-1\right|
\]
 is maximized when 
\[
\Re\left(\frac{1}{z}\right)=\frac{k}{2}.
\]
Thus,
\[
|A_{2}|\leq\sum_{\substack{1\leq k\leq N\\
k\equiv0(\bmod5)
}
}\sum_{\substack{0\leq h<k\\
\gcd(h,k)=1\\
h\equiv1,4(\bmod5)
}
}e^{2\pi\varrho(n-1)}e^{\pi}\frac{2}{KN}\prod_{j=1}^{2}\left|\psi^{\delta_{j}}(\tilde{\varsigma_{j}};\tilde{\tau_{j}})-1\right|
\]

Since for $j=1,2$, 
\[
\lambda_{m_{j},r_{j}}^{*}(\aleph,l)\geq\frac{1}{5},\quad\Im(\tilde{\varsigma_{j}})=\frac{5}{k}\Re\left(\frac{1}{z}\right)\lambda_{m_{j},r_{j}}^{*}(\aleph,l)\geq\frac{5}{2}\lambda_{m_{j},r_{j}}^{*}(\aleph,l)\geq\frac{1}{2}
\]
and 
\[
\Im(\tilde{\tau_{j}})=\frac{5}{k}\Re\left(\frac{1}{z}\right)\geq\frac{5}{2}
\]
we, by \eqref{eq: 3.12}, have

\[
|A_{2}|\leq4\sum_{\substack{1\leq k\leq N\\
k\equiv0(\bmod5)
}
}\sum_{\substack{0\leq h<k\\
\gcd(h,k)=1\\
h\equiv1,4(\bmod5)
}
}e^{2\pi\varrho(n-1)}e^{\pi}\frac{2}{KN}<185e^{2\pi\varrho(n-1)}.
\]

For $A_{1},$ we use Lemma 2.2 to get
\begin{align*}
A_{1} & =\sum_{\substack{1\leq k\leq N\\
k\equiv0(\bmod5)
}
}\sum_{\substack{0\leq h<k\\
\gcd(h,k)=1\\
h\equiv1,4(\bmod5)
}
}e^{-2\pi inh/k}(-1)^{\sum_{j=1}^{J}\delta_{j}\lambda_{m_{j},r_{j}}(h,k)}\omega_{h,k}^{2}\varUpsilon_{h,k}\\
 & \times\Big(\Big(\frac{24n-24}{24}\Big)^{-1/2}I_{-1}\big(\frac{\pi}{6k}\sqrt{24(24n-24)}\big)+E^{*}(I)\Big)
\end{align*}
where 
\[
|E^{*}(I)|<\frac{e^{8\pi}e^{2\pi\varrho(n-1)}}{\pi(n-1)}.
\]
Then, 
\begin{align*}
A(n) & =\sum_{\substack{1\leq k\leq N\\
k\equiv0(\bmod5)
}
}\sum_{\substack{0\leq h<k\\
\gcd(h,k)=1\\
h\equiv1,4(\bmod5)
}
}e^{-2\pi inh/k}(-1)^{\sum_{j=1}^{J}\delta_{j}\lambda_{m_{j},r_{j}}(h,k)}\\
 & \times\omega_{h,k}^{2}\varUpsilon_{h,k}\left((n-1)^{-1/2}I_{-1}\left(\frac{4\pi}{k}\sqrt{n-1}\right)\right)+E^{(1)}(I),
\end{align*}
where 
\begin{align*}
|E^{(1)}(I)| & \leq2e^{2\pi\rho(n-1)+54}+185e^{2\pi\varrho(n-1)}+\sum_{\substack{1\leq k\leq N\\
k\equiv0(\bmod5)
}
}\sum_{\substack{0\leq h<k\\
\gcd(h,k)=1\\
h\equiv1,4(\bmod5)
}
}\frac{e^{8\pi}e^{2\pi\varrho(n-1)}}{\pi(n-1)}\\
 & \leq(2e^{54}+185)e^{2\pi\rho(n-1)}+\frac{N^{2}}{2}\times\frac{e^{8\pi}e^{2\pi\varrho(n-1)}}{\pi(n-1)}\coloneqq A'.
\end{align*}
By Lemma \ref{l2.1}, we have 
\begin{align*}
 & \sum_{\substack{10\leq k\leq N\\
k\equiv0(\bmod5)
}
}\sum_{\substack{0\leq h<k\\
\gcd(h,k)=1\\
h\equiv1,4(\bmod5)
}
}I_{-1}\left(\frac{4\pi}{k}\sqrt{n-1}\right)\\
 & <\sum_{\substack{10\leq k\leq N\\
k\equiv0(\bmod5)
}
}\sum_{\substack{0\leq h<k\\
\gcd(h,k)=1\\
h\equiv1,4(\bmod5)
}
}\sqrt{\frac{\pi}{8}}\frac{e^{\frac{4\pi}{k}\sqrt{n-1}}}{\sqrt{\frac{4\pi}{k}\sqrt{n-1}}}\\
 & <\sum_{\substack{10\leq k\leq N\\
k\equiv0(\bmod5)
}
}\sum_{\substack{0\leq h<k\\
\gcd(h,k)=1\\
h\equiv1,4(\bmod5)
}
}\sqrt{\frac{\pi}{8}}\frac{e^{\frac{4\pi}{10}\sqrt{n-1}}}{\sqrt{\frac{4\pi}{N}\sqrt{n-1}}}\\
 & <\frac{N^{2}}{5}\sqrt{\frac{\pi}{8}}\frac{e^{\frac{2\pi}{5}\sqrt{n-1}}}{\sqrt{\frac{4\pi}{N}\sqrt{n-1}}}\coloneqq B'.
\end{align*}
Hence,
\begin{align*}
\left|\sum_{\substack{10\leq k\leq N\\
k\equiv0(\bmod5)
}
}\sum_{\substack{0\leq h<k\\
\gcd(h,k)=1\\
h\equiv1,4(\bmod5)
}
}e^{-2\pi inh/k}(-1)^{\sum_{j=1}^{J}\delta_{j}\lambda_{m_{j},r_{j}}(h,k)}\omega_{h,k}^{2}\varUpsilon_{h,k}I_{-1}\left(\frac{4\pi}{k}\sqrt{n-1}\right)\right|<B'.
\end{align*}

Let 
\begin{align*}
B(I) & =\sum_{\substack{10\leq k\leq N\\
k\equiv0(\bmod5)
}
}\sum_{\substack{0\leq h<k\\
\gcd(h,k)=1\\
h\equiv1,4(\bmod5)
}
}e^{-2\pi inh/k}(-1)^{\sum_{j=1}^{J}\delta_{j}\lambda_{m_{j},r_{j}}(h,k)}\\
 & \times\omega_{h,k}^{2}\varUpsilon_{h,k}(n-1)^{-\frac{1}{2}}I_{-1}\left(\frac{4\pi}{k}\sqrt{n-1}\right)+E^{(1)}(I).
\end{align*}
Then
\begin{align*}
A(n) & =e^{-2\pi in\times\nicefrac{1}{5}}(-1)^{\sum_{j=1}^{J}\delta_{j}\lambda_{m_{j},r_{j}}(1,5)}\omega_{1,5}^{2}\varUpsilon_{1,5}\\
 & \quad+e^{-2\pi in\times\nicefrac{4}{5}}(-1)^{\sum_{j=1}^{J}\delta_{j}\lambda_{m_{j},r_{j}}(4,5)}\omega_{4,5}^{2}\varUpsilon_{4,5}+B(I)\\
 & =-\frac{4\pi}{5}\cos\left(\frac{\pi}{5}(2n+1)\right)(n-1)^{-\frac{1}{2}}I_{-1}\left(\frac{4\pi}{5}\sqrt{n-1}\right)+B^{(1)}(I),
\end{align*}
where
\[
|B^{(1)}(I)|<|A'|+(n-1)^{-\frac{1}{2}}B'.
\]

Let $N=\lfloor\sqrt{2\pi(n-1)}\rfloor.$ Since $\lfloor\sqrt{2\pi(n-1)}\rfloor>\frac{\sqrt{2\pi(n-1)}}{\sqrt{2}}$
for $n\geq20,$ we have
\begin{align*}
|A'| & \leq(2e^{54}+185)e^{2\pi\rho(n-1)}+\frac{2\pi(n-1)}{2}\times\frac{e^{8\pi}e^{2\pi\varrho(n-1)}}{\pi(n-1)}\\
 & \leq(2e^{54}+e^{8\pi}+185)e^{2}
\end{align*}
and

\[
B'=\frac{N^{\frac{5}{2}}}{20\sqrt{2}}\frac{e^{\frac{2\pi}{5}\sqrt{n-1}}}{(n-1)^{\frac{1}{4}}}<\frac{2\pi^{\frac{5}{4}}e^{\frac{2\pi}{5}\sqrt{n-1}}}{5}(n-1),
\]
where we used the the simple inequality $\lfloor\sqrt{2\pi(n-1)}\rfloor<2\sqrt{\pi(n-1)}.$
Then
\[
|B^{(1)}(I)|<(2e^{54}+e^{8\pi}+185)e^{2}+\frac{2\pi^{\frac{5}{4}}e^{\frac{2\pi}{5}\sqrt{n-1}}}{5}(n-1)^{\frac{1}{2}}
\]
By Lemma \ref{l2.1}, we get 
\[
I_{-1}\left(\frac{4\pi}{5}\sqrt{n-1}\right)>\frac{1}{10}\frac{e^{\frac{4\pi}{5}\sqrt{n-1}}}{\sqrt{\frac{4\pi}{5}\sqrt{n-1}}},\quad n\geq3.
\]
Therefore, when $n>800,$

\begin{align*}
 & \left|-\frac{4\pi}{5}\cos\left(\frac{\pi}{5}(2n+1)\right)(n-1)^{-\frac{1}{2}}I_{-1}\left(\frac{4\pi}{5}\sqrt{n-1}\right)\right|\\
 & >\frac{1}{10}\left|\frac{4\pi}{5}\cos\left(\frac{\pi}{5}(2n+1)\right)\right|(n-1)^{-1/2}\frac{e^{\frac{4\pi}{5}\sqrt{n-1}}}{\sqrt{\frac{4\pi}{5}\sqrt{n-1}}}\\
 & >(2e^{54}+e^{8\pi}+185)e^{2}+\frac{2\pi^{\frac{5}{4}}e^{\frac{2\pi}{5}\sqrt{n-1}}}{5}(n-1)^{\frac{1}{2}},
\end{align*}
so that for $n>160,$ 
\[
A(5n)=-\frac{4\pi}{5}\cos\frac{\pi}{5}\times(5n-1)^{-\frac{1}{2}}I_{-1}\left(\frac{4\pi}{5}\sqrt{5n-1}\right)+B^{(1)}(I)<0.
\]
Calculating the first $800$ coefficients of $1/R^{5}(q)$ by usingthe
software \textit{Maple }we find that $A(5n)<0$ holds for $1\leq n\leq160.$
This proves $A(5n)<0$ for $n\geq1.$

\subsection{Proof of $B(5n)<0$}

Note that

\[
R(q)^{5}=\psi(2\tau;5\tau)^{-5}\psi(\tau;5\tau)^{5}.
\]
Then
\[
m_{1}=m_{2}=5,~r_{1}=2~r_{2}=1,\delta_{1}=-5~\delta_{2}=5.
\]
Hence $L=5,\;\Omega=24$ and $\mathcal{L}_{>0}=\left\{ (2,5),(3,5)\right\} .$
Similarly, we will spilt the integral into two parts according to
$\Delta(\aleph,l)\leq0$ or $\Delta(\aleph,l)>0.$

For $\Delta(\aleph,l)\leq0$, from
\begin{equation}
\left|\int_{\xi_{h,k}}\exp\left(\frac{\pi}{12k}\left(\Omega z+\Delta(\aleph,l)z^{-1}\right)\right)d\varphi\right|<e^{2\pi n\rho}\frac{2}{KN}\label{eq;zhd}
\end{equation}
we know this part of the integral is bounded above, and its upper
bound is $2e^{2\pi\rho(n+1)+54}.$

For $\Delta(\aleph,l)>0,$ proceeding as we handled $S_{\aleph,l}^{(2)},$
we spilt the integral into two parts. For the part of the integrand
involving $\prod_{j=1}^{2}\psi^{\delta_{j}}(\tilde{\varsigma_{j}};\tilde{\tau_{j}})-1,$
the inequality 
\[
5\sum_{j=1}^{2}\left(\lambda_{m_{j},r_{j}}^{*}(\aleph,5)s_{j}+(1-\lambda_{m_{j},r_{j}}^{*}(\aleph,5))t_{j}+n_{j}\right)\geq1
\]
 still holds. This implies that 
\[
\exp\left(\frac{2\pi}{k}\Re\left(\frac{1}{z}\right)\right)\left|\prod_{j=1}^{2}\psi^{\delta_{j}}(\tilde{\varsigma_{j}};\tilde{\tau_{j}})-1\right|
\]
 is maximized when $\Re\left(\frac{1}{z}\right)=\frac{k}{2}.$ Then
the part of this integral can be bounded above and its upper bound
is $185e^{2\pi\varrho(n+1)}.$

The second part is the main term, which can be obtained from Lemma
\ref{12.2}.

From Lemma \ref{12.2}, we deduce that 
\begin{align*}
B(n) & =\sum_{\substack{1\leq k\leq N\\
k\equiv0(\bmod5)
}
}\sum_{\substack{0\leq h<k\\
\gcd(h,k)=1\\
h\equiv2,3(\bmod5)
}
}e^{-2\pi inh/k}(-1)^{\sum_{j=1}^{J}\delta_{j}\lambda_{m_{j},r_{j}}(h,k)}\omega_{h,k}^{2}\varUpsilon_{h,k}\\
 & \quad\times\left((n+1)I_{-1}\left(\frac{4\pi}{k}\sqrt{n+1}\right)\right)+E^{(2)}(I)\\
 & =-\frac{4\pi}{5}\cos\left(\frac{2\pi}{5}(2n-1)\right)(n+1)^{-\frac{1}{2}}I_{-1}\left(\frac{4\pi}{5}\sqrt{n+1}\right)+B^{(2)}(I),
\end{align*}
where 
\[
|E^{(2)}(I)|<(2e^{54}+185)e^{2\pi\rho(n+1)}+\frac{N^{2}}{2}\times\frac{e^{8\pi}e^{2\pi\varrho(n+1)}}{\pi(n+1)}
\]
and 
\[
|B^{(2)}(I)|<|E^{(2)}(I)|+\frac{N^{2}}{5}\sqrt{\frac{\pi}{8}}\frac{e^{\frac{2\pi}{5}\sqrt{n+1}}}{\sqrt{\frac{4\pi}{N}\sqrt{n+1}}}(n+1)^{\frac{1}{2}}.
\]
Setting $N=\lfloor\sqrt{2\pi(n+1)}\rfloor,$we have, for $n\geq20,$
\[
|B^{(2)}(I)|<(2e^{54}+e^{8\pi}+185)e^{2}+\frac{2\pi^{\frac{5}{4}}e^{\frac{2\pi}{5}\sqrt{n+1}}}{5}(n+1)^{\frac{1}{2}}.
\]
It follows from Lemma \ref{l2.1} that 
\[
I_{-1}\left(\frac{4\pi}{5}\sqrt{n+1}\right)>\frac{1}{10}\frac{e^{\frac{4\pi}{5}\sqrt{n+1}}}{\sqrt{\frac{4\pi}{5}\sqrt{n+1}}},\quad n\geq1.
\]
Then, when $n>800,$
\begin{align*}
 & \left|-\frac{4\pi}{5}\cos\left(\frac{2\pi}{5}(2n-1)\right)(n+1)^{-\frac{1}{2}}I_{-1}\left(\frac{4\pi}{5}\sqrt{n+1}\right)\right|\\
 & >\left|\frac{4\pi}{5}\cos\left(\frac{2\pi}{5}(2n-1)\right)(n+1)^{-\frac{1}{2}}\times\frac{1}{10}\frac{e^{\frac{4\pi}{5}\sqrt{n+1}}}{\sqrt{\frac{4\pi}{5}\sqrt{n+1}}}\right|\\
 & >(2e^{54}+e^{8\pi}+185)e^{2}+\frac{2\pi^{\frac{5}{4}}e^{\frac{2\pi}{5}\sqrt{n+1}}}{5}(n+1)^{\frac{1}{2}}.
\end{align*}
Therefore, when $n>160,$ 
\[
B(5n)=-\frac{4\pi}{5}\cos\frac{2\pi}{5}(5n+1)^{-\frac{1}{2}}I_{-1}\left(\frac{4\pi}{5}\sqrt{5n+1}\right)+B(I)<0.
\]
Using the software \textit{Maple }to calculate the first $800$ coefficients
of $R^{5}(q),$ we find that $B(5n)<0$ is true for $1\leq n\leq160.$
This shows $B(5n)<0$ for $n\geq1.$

\subsection{Proof of $D(5n+1)>0$}

Note that 
\[
\frac{R(q^{5})}{R(q)^{5}}=\psi(2\tau;5\tau)^{5}\psi(\tau;5\tau)^{-5}\psi(5\tau;25\tau)\psi(10\tau;25\tau)^{-1},
\]
Then,
\[
m_{1}=m_{2}=5,\ m_{3}=m_{4}=25,r_{1}=2,\ r_{2}=1,\ r_{3}=5,~r_{4}=10,
\]
 and
\[
\delta_{1}=-5,~\delta_{2}=5,\ \delta_{3}=1,~\delta_{4}=-1.
\]
Hence $L=25,\;\Omega=0,$ and $\Delta(\aleph,l)>0$ occurs only when
$\aleph\equiv1,4(\bmod5)$ and $l\equiv5,10,15,20(\bmod25).$ In particular,
$\Delta(\aleph,l)=24$ when $\Delta(\aleph,l)>0.$ 

For $\Delta(\aleph,l)\leq0$, by \ref{eq;zhd}, we can determine the
upper bound of this part of the integral and its upper bound is $e^{2\pi n\rho+330}.$

For $\Delta(\aleph,l)>0,$ we spilt the integral into two parts. For
the part of the integrand involving $\prod_{j=1}^{4}\psi^{\delta_{j}}(\tilde{\varsigma_{j}};\tilde{\tau_{j}})-\Pi_{h,k},$
we know that $\Pi_{h,k}\neq1$, which differs from the previous ones
in the proofs of $A(5n)<0$ and $B(5n)<0.$

By \eqref{eq: 3.10}, we get

\begin{align*}
\prod_{j=1}^{4}\psi^{\delta_{j}}(\tilde{\varsigma_{j}};\tilde{\tau_{j}}) & =\Pi_{h,k}\left(\prod_{j=1}^{2}(e^{2\pi i\tilde{\varsigma_{j}}},e^{2\pi i(\tilde{\tau_{j}}-\tilde{\varsigma_{j}})};e^{2\pi i\tilde{\tau_{j}}})_{\infty}^{\delta_{j}}\right.\\
 & \times\left.\prod_{j=3}^{4}(e^{2\pi i\tilde{(\varsigma_{j}}+\tilde{\tau_{j}})},e^{2\pi i(\tilde{\tau_{j}}-\tilde{\varsigma_{j}})};e^{2\pi i\tilde{\tau_{j}}})_{\infty}^{\delta_{j}}\right)
\end{align*}

Let 
\[
\zeta_{j}^{*}:=\begin{cases}
\tilde{\varsigma_{j}} & \text{if }j=1,2\\
\tilde{\varsigma_{j}}+\tilde{\tau_{j}} & \text{if }j=3.4.
\end{cases}
\]
Then

\[
\Im(\zeta_{j}^{*})=\Phi(\lambda_{m_{j},r_{j}}^{*}(h,k))\frac{\gcd^{2}(m_{j},\ell)}{m_{j}k}\Re(z^{-1}),
\]
where for $0\leq x<1,$

\[
\Phi(x):=\begin{cases}
1, & \text{if }x=0,\\
x, & \text{otherwise}.
\end{cases}
\]
Notice that 
\begin{align*}
 & \left|\prod_{j=1}^{4}\psi^{\delta_{j}}(\tilde{\varsigma}_{j};\tilde{\tau}_{j})-\prod_{h,k}\right|\\
 & =\prod_{j=1}^{4}\left|\left(e^{2\pi i\tilde{\varsigma}_{j}},e^{2\pi i(\tilde{\tau}_{j}-\tilde{\varsigma}_{j})};e^{2\pi i\tilde{\tau}_{j}}\right)_{\infty}^{\delta_{j}}-\prod_{h,k}\right|\\
 & =\prod_{j=1}^{4}\left|\prod_{h,k}\right|\left|\frac{1}{\prod_{h,k}}\left(e^{2\pi i\tilde{\varsigma}_{j}},e^{2\pi i(\tilde{\tau}_{j}-\tilde{\varsigma}_{j})};e^{2\pi i\tilde{\tau}_{j}}\right)_{\infty}^{\delta_{j}}-1\right|\\
 & =\left|\prod_{h,k}\right|\left|\prod_{j=1}^{4}\left(e^{2\pi i\zeta_{j}^{*}},e^{2\pi i(\tilde{\tau}_{j}-\tilde{\varsigma}_{j})};e^{2\pi i\tilde{\tau}_{j}}\right)_{\infty}^{\delta_{j}}-1\right|.
\end{align*}
 Since
\[
\sum_{j=1}^{4}\frac{\gcd^{2}(m_{j},l)}{m_{j}}\left(\Phi(\lambda_{m_{j},r_{j}}^{*}(h,k))s_{j}+\left(1-\lambda_{m_{j},r_{j}}^{*}(\aleph,l)\right)t_{j}+n_{j}\right)\geq\frac{\Delta(\aleph,l)}{24}
\]
we know that

\[
\exp\left(\frac{\Delta(\aleph,l)\pi}{12k}\Re\left(\frac{1}{z}\right)\right)\left|\prod_{j=1}^{4}\psi^{\delta_{j}}(\tilde{\varsigma}_{j};\tilde{\tau}_{j})-\prod_{h,k}\right|
\]
is maximized when 
\[
\Re\left(\frac{1}{z}\right)=\frac{k}{2},
\]
and an upper bound of this part of the integral is $e^{2\pi n\rho+270}.$

For the other part, by Lemmas \ref{12.2} and \ref{l2.1}, we know
that this part produces the dominant term.

Form Lemmas \ref{12.2} and \ref{l2.1}, we derive that 
\begin{align*}
D(n) & =e^{-2\pi in\nicefrac{1}{5}}(-1)^{\sum_{j=1}^{4}\delta_{j}\lambda_{m_{j},r_{j}}(1,5)}\omega_{1,5}^{2}\varUpsilon_{1,5}\Pi_{1,5}\times n{}^{-\frac{1}{2}}I_{-1}\left(\frac{4\pi}{5}\sqrt{n}\right)\\
 & +e^{-2\pi in\nicefrac{4}{5}}(-1)^{\sum_{j=1}^{4}\delta_{j}\lambda_{m_{j},r_{j}}(4,5)}\omega_{4,5}^{2}\varUpsilon_{4,5}\Pi_{4,5}\times n^{-\frac{1}{2}}I_{-1}\left(\frac{4\pi}{5}\sqrt{n}\right)+B^{(3)}(I)\\
 & =-\frac{2\pi}{5}\frac{\cos\frac{\pi}{5}}{1+\cos\frac{2\pi}{5}}n{}^{-\frac{1}{2}}\cos\left(\frac{\pi}{5}(2n+1)\right)I_{-1}\left(\frac{4\pi}{5}\sqrt{n}\right)+B^{(3)}(I),
\end{align*}
 Where 
\[
|B^{(3)}(I)|<e^{2\pi n\rho+330}+e^{2\pi n\rho+270}+\frac{N^{2}}{2}\times\frac{e^{8\pi}e^{2\pi n\varrho}}{n\pi}+\frac{N^{2}}{5}\sqrt{\frac{\pi}{8}}\frac{e^{\frac{2\pi}{5}\sqrt{n}}n^{\frac{1}{2}}}{\sqrt{\frac{4\pi}{N}\sqrt{n}}}.
\]

\noindent Taking $N=\lfloor\sqrt{2\pi n}\rfloor,$ we have, for $n\geq20,$

\[
|B^{(3)}(I)|<e^{332}+e^{272}+e^{8\pi+2}+\frac{2\pi^{\frac{5}{4}}e^{\frac{2\pi}{5}\sqrt{n}}n^{\frac{1}{2}}}{5}.
\]
It is deduced from Lemma \ref{l2.1} that 
\[
I_{-1}\left(\frac{4\pi}{5}\sqrt{n}\right)>\frac{1}{10}\frac{e^{\frac{4\pi}{5}\sqrt{n}}}{\sqrt{\frac{4\pi}{5}\sqrt{n}}},\quad n\geq2.
\]
Then, when $n>19000,$ 
\begin{align*}
 & \left|-\frac{2\pi}{5}\frac{\cos\frac{\pi}{5}}{1+\cos\frac{2\pi}{5}}n{}^{-\frac{1}{2}}\cos\left(\frac{\pi}{5}(2n+1)\right)I_{-1}\left(\frac{4\pi}{5}\sqrt{n}\right)\right|\\
 & >\frac{2\pi}{5}\frac{\cos\frac{\pi}{5}}{1+\cos\frac{2\pi}{5}}n{}^{-\frac{1}{2}}\left|\cos\left(\frac{\pi}{5}(2n+1)\right)\right|\times\frac{1}{10}\frac{e^{\frac{4\pi}{5}\sqrt{n}}}{\sqrt{\frac{4\pi}{5}\sqrt{n}}}\\
 & >e^{332}+e^{272}+e^{8\pi+2}+\frac{2\pi^{\frac{5}{4}}e^{\frac{2\pi}{5}\sqrt{n}}n^{\frac{1}{2}}}{5}.
\end{align*}
Therefore, when $n\geq3800,$
\[
D(5n+1)=-\frac{2\pi}{5}\frac{\cos\frac{\pi}{5}\cos\frac{3\pi}{5}}{1+\cos\frac{2\pi}{5}}(5n+1)^{-\frac{1}{2}}I_{-1}\left(\frac{4\pi}{5}\sqrt{5n+1}\right)+B^{(3)}(I)>0.
\]
We employ the software \textit{Maple }to calculate the first $19000$
coefficients of $R(q^{5})/R^{5}(q)$ and find that $D(5n+1)>0$ holds
for $0\leq n\leq3800.$ This proves $D(5n+1)>0$ for all $n\geq0.$

\end{document}